\documentclass{amsart}
\usepackage{euscript,amsmath, amssymb, amsfonts}

\numberwithin{equation}{section}

\usepackage[left=2cm, right=2cm]{geometry}
\usepackage{cite}
\usepackage{graphicx}
\usepackage{hyperref}
\usepackage{comment}
\usepackage[numbers,square]{natbib}

\theoremstyle{definition}

\begin{document}

\title{Generating functions and efficient numerical simulation of aggregation-fragmentation kinetics}

\author{S.~A.~Matveev}
\address{Faculty of Computational Mathematics and Cybernetics, Lomonosov Moscow State University, Leninskiye gory GSP-1, Moscow, 119991,  Russia; 
Marchuk Institute of Numerical Mathematics, Russian Academy of Science, Gubkin st.~8, Moscow, 119333, Russia}
\email{matseralex@cs.msu.ru}

\keywords{aggregation, fragmentation, exact solutions, kinetic equations, generating function, numerical simulation, low-rank matrices}
\subjclass[2020]{34A34, 35B40, 82C40}

\begin{abstract}
In this work, we review and revisit the generating function techniques that provide exact analytical solutions for aggregation and fragmentation equations across several physical regimes including spontaneous and collisional shattering. For discrete coagulation-fragmentation equations with size-independent rates under monodisperse initial conditions, we show the derivation of several explicit closed-form solutions. We also briefly report the exact solutions for continuous, three-particle, $D$-particle collisions and two-component generalizations. Source-driven aggregation yields steady distributions featuring a universal $s^{-3/2}$ power-law decay and a cutoff mass scaling $s_{*} \sim t^{2}$. We also remind that exploitation of the low-rank representations of the aggregation kernel allows to organize efficient modelling the aggregation kinetics and extend this idea showing a simple way for the fast accounting of the spontaneous framentation terms. We exploit our numerical algorithm and study the numerical solutions for multiple models of spontaneous fragmentation and find a new example of temporal oscillations for the closed aggregation-fragmentation process. 
\end{abstract}

\maketitle

\section{Introduction}\label{sec1}

The kinetics of aggregation and fragmentation processes often can be described using systems of nonlinear differential equations \cite{smoluchowski1916drei, Galkin2001, Krapivsky}. The smallest basic particles are usually called monomers, and their size is formally taken as unity. Collisions of small particles can form larger aggregates, which contain several monomers at once. Further collisions of larger aggregates can also form increasingly larger particles, consisting of thousands or even millions of monomers. In the case where the aggregation coefficients (cores) $K_{i,j}$ for a reaction of the type $[i]+[j] \rightarrow [i+j]$ are known or specified using a formal expression (for example, $K_{i,j}=(i/j)^{1/3}+(j/i)^{1/3}+2$ \cite{Leyvraz}), the system of equations describing the process of aggregate growth is well known and is called the Smoluchowski equations \cite{smoluchowski1916drei}:
\begin{equation*}\label{eq:Smolbasic}
\frac{d n_s}{dt} = \underbrace{\frac{1}{2} \sum_{i+j=s} K_{i,j} n_i n_j}_{\text{formation of particles of size } s} - \underbrace{n_s \sum_{i =1}^{\infty} K_{s,i} n_i,}_{\text{disappearance of particles of size } s}s = 1,2,\ldots,\infty.
\end{equation*}
These equations describe the time dynamics of the concentration functions $n_s(t)$ of particles of size $s$ due to their formation in collisions with smaller particles and their disappearance in collisions with other particles. Within this model, all collisions are assumed to be paired, satisfying the basic law of conservation of mass. It is additionally assumed that the particles are spatially uniformly distributed \cite{smoluchowski1916drei, Galkin2001, Krapivsky, Leyvraz}.

For given initial conditions on the concentrations $n_s(t=0)$, the Cauchy problem arises for a formally infinite system of nonlinear Smoluchowski kinetic equations. For a sufficiently wide class of symmetric homogeneous aggregation kernels
$$
0 \leq K_{\ell i,\ell j} = K_{\ell j,\ell i} = \ell^{\rho} K_{i,j}
$$
with homogeneity index $\rho \leq 1$, it is proved that the solution to the Cauchy problem exists, is unique, and continuously depends on the input data \cite{Galkin2001, Melzak1957, mcleod62}. For a wide class of non-negative initial conditions (for more details, see Chapter 2 of \cite{Galkin2001}), the solution to the Cauchy problem will preserve non-negativity, will be uniformly bounded, and will satisfy the property of conservation of the total mass of matter.
$$
\sum_{s=1}^{\infty} s n_s(t) = \sum_{s=1}^{\infty} s n_s(0) \equiv \text{CONST}.
$$
This class includes the following kernels, which are important for analytical research and various applications:
\begin{enumerate}
\item Constant $K_{i,j} = K$ \cite{smoluchowski1916drei},
\item Additive $K_{i,j} = i+j$ \cite{lushnikov1978coagulation},
\item Generalized product $K_{i,j}= (i \cdot j)^{\mu}$, $\mu < 1/2$ \cite{hayakawa1987irreversible, Colm},
\item Ballistic $K_{i,j}= (i^{1/3} + j^{1/3}) \sqrt{\frac{1}{i}+\frac{1}{j}}$ \cite{brilliantov2009model},
\item Diffusion $K_{i,j}= (i/j)^{1/3} + (j/i)^{1/3} + 2$ \cite{krapivsky2012driven},
\item $\alpha-\beta$ model $K_{i,j} = i^{\alpha} j^{\beta} + i^{\beta} j^{\alpha}$ for $\alpha+\beta < 1$, $\alpha <1$, $\beta < 1$ \cite{Colm, hayakawa1987irreversible}.
\end{enumerate}
For monodisperse initial conditions $n_s(t=0) = \delta_{s,1}$ (where $\delta_{s,1}$ is Kronecker symbol) exact analytical solutions are known \cite{smoluchowski1916drei, lushnikov1978coagulation, Lushnikov04} for the Cauchy problems with kernels 1 -- 3 and in this work we remind a simple way to obtain them via the generating function method for the sequence $n_s(t)$.

The main aim of this paper is to review how generating function methods yield exact solutions for coagulation-fragmentation models, including source-driven and multi-particle systems, and to characterize universal asymptotic behaviors such as the $s^{-3/2}$ power-law decay. 

We also verify our analytical observations numerically. We discuss an efficient way for organizing the efficient computations in the case with spontaneous fragmentation terms and study structure of the stationary solutions for several models. All in all, we find a new interesting example of non-decaying temporal oscillations for aggregation-fragmentaion equations.

\section{Basic irreversible coagulation}

\subsection{Constant coagulation kernel}

Here we concentrate on a problem with a constant aggregation kernel, e.g., $K_{i,j} = 2$, we introduce the generating function for the sequence $n_s(t)$
$$
g(z,t) = \sum_{s=1}^{\infty} n_s(t) z^s.
$$
Before working with the generating function, we first derive an elementary differential equation for the total aggregate density $n(t) = \sum\limits_{s=1}^{\infty} n_s(t)$
\begin{equation*}
\frac{dn(t)}{dt} = \sum_{s=1}^{\infty} \frac{dn_s(t)}{dt} = -n^2, \qquad n(t=0)=1.
\end{equation*}
From this equation, it follows that $n(t) = (1+t)^{-1}$. Next, we note that
\begin{align*}
& \frac{\partial g(z,t)}{\partial t} = \sum_{s=1}^{\infty} \frac{dn_s}{dt} z^s \nonumber
= \sum_{s=2}^{\infty} \left( \sum_{i=1}^{s-1} n_i z^i n_{s-i} z^{s-i} \right) - 2n(t) \sum_{s=1}^{\infty} n_s z^s \nonumber\\
&= g^2(z,t) - 2n(t)g(z,t) = g(z,t)\bigl(g(z,t) - n(t)\bigr), \qquad g(z,t=0)=z.
\end{align*}
Further, we make the substitution $\hat{g}(z,t) = g(z,t) - n(t)$ and get a parametric equation
\begin{equation*}
\frac{\partial\hat{g}(z,t)}{\partial t} = \hat{g}^2(z,t), \qquad \hat{g}(z,t=0) = z-1.
\end{equation*}
Using separation of variables, we obtain the solution $\hat{g}(z,t) = \dfrac{z-1}{1 - t(z-1)}$. Therefore,
\begin{equation*}
g(z,t) = \hat{g}(z,t) + n(t) = \frac{z-1}{1 - t(z-1)} + \frac{1}{1+t} = \frac{z}{1 - t(z-1)} \cdot \frac{1}{1+t}. 
\end{equation*}
Expansion of $\dfrac{z}{1 - t(z-1)}$ as a Taylor series in powers of $z$ gives
\begin{equation*}
g(z,t) = \frac{1}{1+t} \sum_{s=1}^{\infty} \frac{t^{s-1}}{(1+t)^s} z^s = \sum_{s=1}^{\infty} \underbrace{\frac{t^{s-1}}{(1+t)^{s+1}}}_{n_s(t)} z^s.
\end{equation*}
Thus, for the Cauchy problem with a constant kernel $K_{i,j}=2$ and monodisperse initial conditions, the analytical expressions for $n_s(t)$ are
\begin{equation*}
n_s(t) = \frac{t^{s-1}}{(1+t)^{s+1}}.
\end{equation*}
Using the homogeneity of the right-hand side, we obtain the solution for an arbitrary constant kernel $K_{i,j}=K$:
\begin{equation*}
n_s(t) = \frac{\left(\frac{K}{2}t\right)^{s-1}}{\left(1 + \frac{K}{2}t\right)^{s+1}}.
\end{equation*}
One may find a transformation (see, e.g., the book \cite{Krapivsky}) that allows one to obtain an exact solution to the Cauchy problem with arbitrary initial conditions based on the analytical solution for monodisperse initial conditions. Despite the existence of such a formal transformation, exact analytical formulas even for the constant kernel can only be derived for certain special classes of initial conditions.

\subsection{More complex cases}

A suitable choice of generating functions for the sequence $n_s(t)$ allows obtaining the analytical solutions similarly (though with somewhat more technical effort) for the additive kernel $K_{i,j}=i+j$ (using the same generating function) and for the multiplicative kernel $K_{i,j}=i\cdot j$ (which requires an exponential generating function).

For the additive kernel with monodisperse initial conditions, using the same generating function $g(z,t) = \sum\limits_{s=1}^{\infty} n_s(t) z^s$, we have the following equations (see, e.g., \cite{fernandez2007exact, lepek2019exact, lepek2026mini}):
\begin{align*}
\frac{\partial g(z,t)}{\partial t} &= z\bigl(g(z,t) - n(t)\bigr)\frac{\partial g}{\partial z} - g(z,t), \qquad \left. z\frac{\partial g}{\partial z}\right|_{z=1} = 1,\\
\frac{dn(t)}{dt} &= -n(t), \qquad n(0)=1,\\
\frac{\partial\hat{g}(z,t)}{\partial t} &= z \hat{g}(z,t) \frac{\partial\hat{g}(z,t)}{\partial z} - \hat{g}(z,t), \quad \text{where } \hat{g}(z,t) = g(z,t) - n(t). 
\end{align*}
Clearly, $n(t) = e^{-t}$. 

At first, the equation for $\hat{g}(z,t)$ can be solved using the method of characteristics and yield the implicit equation for $\hat{g}(z,t)$
$$
z e^{\hat{g}(z,t) (e^{t} - 1)} = \hat{g}(z,t) e^{t} + 1.
$$
Then, we assume that $g(z,t) = n(t) \cdot \tau(z,t)$ and see that $\tau(z,t)$ satisfies the transcendental equation
\begin{equation*}
z = \tau(z,t)\, e^{\alpha(t)(1 - \tau(z,t))}, \qquad \alpha(t) = 1 - n(t).
\end{equation*}
One may note that $\tau(z,t)$ is the ratio $g(z,t)/n(t)$ and $\alpha(t) = 1 - n(t) $ is a known function of time. Hence, we need to expand $\tau(z,t)$ as a power series in $z$
$$
\tau(z,t) = \sum_{s=1}^{\infty} a_s(\alpha(t) ) z^{s}.
$$
The equation $z = \tau e^{\alpha(t) \cdot (1 - \tau(z,t))}$ can be rewritten in a suitable form for inversion. We multiply both sides of equation by $-\alpha(t) \cdot e^{-\alpha(t)}$ and set $w = -\alpha(t) \tau$. Then we obtain 
$$w e^{w} = -\alpha z e^{-\alpha(t)}$$
and utilize the Lambert $W$-function
$
w = W\!\left(-\alpha(t) z e^{-\alpha(t)}\right).
$
Hence, we obtain 
$$
 \tau(z,t) = -\frac{1}{\alpha(t)} W\!\left(-\alpha(t) z e^{-\alpha(t) }\right).
$$
Further, we use the Lagrange inversion theorem directly for expansion of series \cite{lavrentyev1969intro, surya2023lagrange} and see
$$
\tau(z,t) = \sum_{m=1}^{\infty} \frac{(m\alpha(t) )^{m-1}}{m!} e^{-m\alpha(t)} z^{m}.
$$
Then, we recall that $g(z,t) = n(t) \tau(z,t)$ and $n(t) = e^{-t}$. Therefore,
$$
g(z,t) = e^{-t} \sum_{m=1}^{\infty} \frac{(m\alpha(t))^{m-1}}{m!} e^{-m\alpha(t)} z^{m}.
$$
Hence,
$$
n_{s}(t) = e^{-t} \cdot \frac{(s\alpha)^{s-1}}{s!} e^{-s\alpha}.
$$
Since $\alpha(t) = 1 - e^{-t}$, we obtain the final well-known analytical solution
$$
n_{s}(t) = e^{-t} \cdot \frac{\big(s(1 - e^{-t})\big)^{s-1}}{s!} e^{s(1 - e^{-t})}.
$$

\begin{figure}
\includegraphics[width=0.65\textwidth]{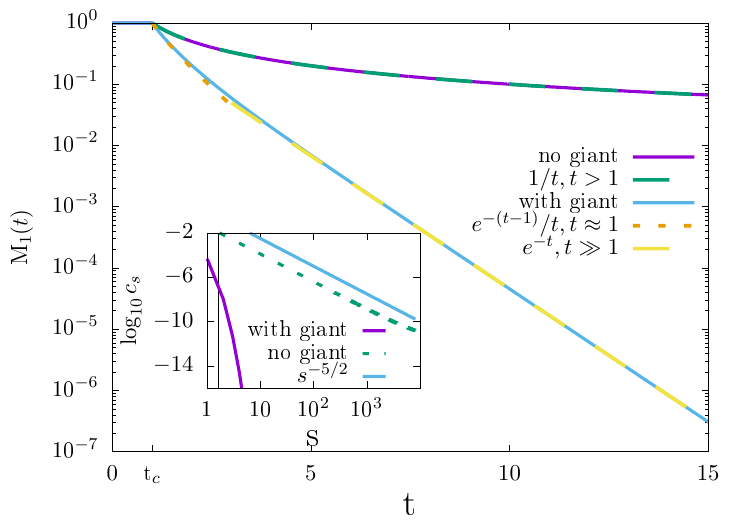}
\caption{Dynamics of the total mass and structure of the particle size distribution at $t=5$ depending on the problem statement with product kernel.}
\label{Fig:gelation}
\end{figure}

In a number of applications, the product kernel $K_{i,j} = i\cdot j$, which is not included in the mass-conserving class (its homogeneity factor $\rho = 2$), plays an important role \cite{lushnikov2006gelation, riordan2016convergence}. For this kernel, a sol-gel transition occurs, meaning that the conservation of total mass is violated. Instead of the ordinary generating function, it is convenient to use an exponential generating function
\begin{equation*}
\mathcal{E}(z,t) = \sum_{s\geq 1} s n_s(t) e^{sz} 
\end{equation*}
associated with the sequence $s n_s(t)$. For the Cauchy problem with monodisperse initial conditions, an analytical solution is also known for this kernel
$$
n_s(t) = \left\lbrace
\begin{matrix}
\dfrac{s^{s-2}}{s!} t^{\,s-1} e^{-s t}, & 0 \le t \le 1, & \text{pre-gel phase} \\
\dfrac{s^{s-2}}{s!} \dfrac{1}{t} \, e^{-s}, & t > 1 & \text{post-gel phase},\\
\end{matrix} \right.
\qquad s = 1, 2, 3, \dots
$$
and can be expressed with use of this generating function \cite{lushnikov2006gelation, krapivsky2024gelation, lepek2019exact}. The total mass of particles decays exponentially in the post-gel phase.

It is worth noting that two types of setting the kinetic equations are available for the product kernel. One may take into account the mass conservation law $\sum_{s=1}^{\infty} s n_s(t) =1$ and analyze the equations as
\begin{eqnarray*}
\frac{dn_s}{dt} = \frac{1}{2} \sum_{i+j=s} i j n_i n_j - s n_s, \qquad s = 1, 2, \ldots, \text{ with active giant component}.
\end{eqnarray*}
The analytical expressions from above are correct for these equations and numerical integration of their finite subsystems allows to capture and approximate the solution.
One may also think about the sequence of the truncated systems assuming $N \rightarrow \infty$ and solve numerically
\begin{eqnarray*}
\frac{dn_s}{dt} = \frac{1}{2} \sum_{i+j=s} i j n_i n_j - s n_s \sum_{i=1}^{N} i n_i, \qquad  s = 1, 2, \ldots, \text{ no active giant component} .
\end{eqnarray*}
Both types of systems correspond to the similar solutions before the gelation time. Surprisingly, the total mass of the sequence of the solutions of such truncated systems decays as $1/t$ in the post-gel phase (see Fig. \ref{Fig:gelation}, ). Such a difference between two formulations of kinetic equations is not possible for the mass-conserving kernels.

This issue becomes even more drastic for the source-enhanced aggregation kinetics with the product coefficients and its multi-particle generalizations \cite{krapivsky1991aggregation}. In the first case the gel engulfes the system and the total mass decays as $1/t$. In the second case, there exist the stationary-distributions with power-law asymptotics $n_s \simeq s^{-5/2}$. Recently, some similar kind of the phase-transition has been analyzed for the multi-particle kernels corresponding to Achlioptas processes \cite{achlioptas2009explosive, riordan2011explosive, riordan2025phase} in context of explosive percolation in random networks.

\subsection{Generalizations}

The classical Smoluchowski equations have been generalized \cite{Galkin2001, muller1928allgemeinen} to the case of a continuous particle size distribution $n(v,t)$, $v>0$
\begin{equation*}
\frac{\partial n(v,t)}{\partial t} = \frac{1}{2} \int_0^v K(v-u,u)\, n(v-u,t)\, n(u,t)\, du - n(v,t) \int_0^\infty K(v,u)\, n(u,t)\, du, 
\end{equation*}
where the integral operators on the right-hand side similarly represent the balance between the formation and depletion of particles of size $v$. The mass conservation law  takes the form
\begin{equation*}
\int_0^\infty v\, n(v,t)\, dv = \int_0^\infty v\, n(v,0)\, dv = \text{CONST}.
\end{equation*}
For the same kernels (constant, additive, and multiplicative) the analytical solutions to this continuous equation can be obtained for exponential initial conditions $n(v,0) = e^{-v}$ using the Laplace transform \cite{Galkin2001, Lushnikov1976, fernandez2007exact}.

For example, the exact solution of the Cauchy problem for the continuous Smoluchowski equation with a constant kernel $K(u,v)=1$ and an exponential initial condition has a structure similar to its discrete counterpart:
\begin{equation*}
n(v,t) = \frac{e^{-v/(1+t/2)}}{(1+t/2)^2}.
\end{equation*}

The solvability theory for the Cauchy problem with bounded kernels and kernels having homogeneity index $\rho < 1$ was first established by Melzak \cite{Melzak1957} and McLeod \cite{mcleod62, mcleod64}. For such kernels, the Cauchy problem for the Smoluchowski equation is well-posed: the solution exists, is unique, and depends continuously on the initial data. Further refinements and a systematic exposition of this theory can be found, for instance, in the monograph by Galkin \cite{Galkin2001}.

For coagulation kernels with higher homogeneity indices ($\rho > 1$), the mass conservation law no longer holds also for the continuous coagulation equation \cite{ball2011instantaneous}. Instead, a sol-gel transition is observed \cite{Flory41, Flory41b, Flory41c} meaning the formation of a gel \cite{Lushnikov04, lushnikov2006gelation} (which, by analogy with growth processes in random graphs, is also commonly called the giant component \cite{ER60, ben2012popularity}), gradually absorbing the entire mass of the dispersed phase \cite{Stockmayer43, Stockmayer44}.

For complex particles consisting of a mixture of several components (e.g., water, sulfuric acid, soot, etc.), the coagulation integro-differential equation is written in multidimensional form \cite{fernandez2007exact, Palaniswaamy2006}:
\begin{equation*}
\begin{aligned}
\frac{\partial n(\overline{v},t)}{\partial t} = \frac{1}{2} \int_0^{v_1} \cdots \int_0^{v_d} K(\overline{v}-\overline{u};\overline{u})\, n(\overline{v}-\overline{u},t)\, n(\overline{u},t)\, du_1\ldots du_d\\
- n(\overline{v},t) \int_0^\infty \cdots \int_0^\infty K(\overline{u};\overline{v})\, n(\overline{u},t)\, du_1\ldots du_d,
\end{aligned}
\end{equation*}
and the law of mass conservation becomes
\begin{equation*}
\int_0^\infty \cdots \int_0^\infty (v_1 + \ldots + v_d)\, n(v_1,\ldots,v_d,t)\, dv_1\ldots dv_d \equiv \text{CONST}.
\end{equation*}

Analytical solutions for multidimensional Smoluchowski equations are also very rare. By analogy with the continuous one-component case, they are typically obtained using the multidimensional Laplace transform \cite{lushnikov1978coagulation, fernandez2007exact}. Several generalizations of standard coagulation kernels are known:
\begin{align*}
K(\overline{u};\overline{v}) &\equiv 1 \quad \text{(constant kernel)},\\
K(\overline{u};\overline{v}) &= \sum_{i=1}^d u_i + \sum_{i=1}^d v_i \quad \text{(generalized additive kernel)},\\
K(\overline{u};\overline{v}) &= \left(\sum_{i=1}^d u_i\right)^{\nu} \left(\sum_{i=1}^d v_i\right)^{\nu} \quad \text{(generalized multiplicative kernel)}.
\end{align*}
The Cauchy problem for the two-component Smoluchowski equation with a constant kernel and initial conditions of the form
\begin{equation*}
n(v_1,v_2,t=0) = ab\, e^{-a v_1 - b v_2}
\end{equation*}
has the following well-known analytical solution (see, e.g., \cite{fernandez2007exact}):
\begin{equation*}
n(v_1,v_2,t) = \frac{ab\, e^{-a v_1 - b v_2}}{(1+t/2)^2}\, I_0\!\left(2\sqrt{\frac{ab\, v_1 v_2\, t}{t+2}}\right),
\end{equation*}
where $a,b$ are arbitrary positive numbers, and $I_0$ is the modified Bessel function of the first kind of order zero. Interestingly, the total particle density in this case can also be calculated analytically and agrees with the much simpler one-component case:
\begin{equation*}
N(t) = \int_0^\infty \int_0^\infty \frac{ab\, e^{-a v_1 - b v_2}}{(1+t/2)^2}\, I_0\!\left(2\sqrt{\frac{ab\, v_1 v_2\, t}{t+2}}\right) dv_1 dv_2 \equiv \frac{1}{1+t/2}.
\end{equation*}

In case of the two-component additive kernel and the same exponential initial conditions the analytical solution can be expressed with rapidly converging series \cite{fernandez2007exact}
$$
n(v_1,v_1,t) = ab(1 - t)\exp \left(-a v_1 - b v_2 - \frac{ab}{a + b} (v_1 + v_2)t\right) \cdot \sum_{k = 0}^{\infty} \dfrac{\left( (a^2 b^2 (v_1 + v_2) v_1 v_2 t) \cdot (a + b)^{-1} \right)^k}{(k + 1)! (k!)^2}.
$$
There are also several other results of analytical solutions expressed as for the non-symmetric additive kernel in the case of the two-component equation \cite{yang2016analytical}, as well as several more complex analytical estimates \cite{piskunov2013analytical, Leyvraz}.

\section{Irreversible source-driven coagulation with constant and slowly growing kernels}
\label{par:Irreversible Coagulation with Sources}


Let us now consider a model of particle aggregation with a constant kernel and without particle sinks:
\begin{equation*}\label{eq:ConstKer}
\frac{d n_s}{dt} = \frac{1}{2} \sum_{i=1}^{s-1} K n_i n_{s-i} - n_s \sum_{i=1}^{\infty} K n_i + P_s,
\end{equation*}
\begin{eqnarray*}
\notag
P_s = \left\lbrace
\begin{matrix}
1, & s = 1, \\
P, & s = \hat{s}, \\
0, & s \neq 1, \hat{s}.
\end{matrix}\right.
\end{eqnarray*}
In this model, the decrease in the concentrations $n_s(t)$ due to particle merging into larger aggregates is compensated by a constant influx of matter from the sources $P_s$ (in general, sources may also depend on time).

If $P = 0$, $K=2$ (only one monomer source is active), then the generating function method readily yields an exact solution to the problem \cite{Krapivsky}. For the dynamics of the total particle concentration, we have
\begin{equation*}
\label{eq:ct-2}
\frac{d n}{dt} = -n^2 + 1,
\end{equation*}
and for a system with zero initial conditions, one easily obtains
\begin{equation*}
\label{eq:ct-sol}
n(t) = \tanh(t).
\end{equation*}
Next, we substitute this result for the total density into the equation for monomers
$$
\frac{d n_1(t)}{dt} = -2\tanh(t) \, n_1(t) + 1,
$$
from which we obtain an analytical expression for the monomer dynamics
\begin{equation*}
\label{c1t-sol}
n_1(t) = \frac{1}{2} \left[ \frac{t}{\cosh^2(t)} + \tanh(t) \right].
\end{equation*}
Using $n(t) = \tanh(t)$, we rewrite the equations for $n_s(t)$ for $s \geq 2$ as
\begin{equation*}
\frac{d n_s}{dt} + 2\tanh(t) \, n_s = \sum_{i+j=s} n_i n_j.
\end{equation*}
After some manipulations of the left-hand side, we obtain
\begin{equation*}
\frac{d}{dt} \left[ n_s(t) \, \cosh^2(t) \right] = \cosh^2(t) \sum_{i+j=s} n_i(t) n_j(t),
\end{equation*}
which, upon integration, yields the relation
\begin{equation*}
\label{ckt-sol}
n_s(t) = \frac{1}{\cosh^2(t)} \int_0^t \cosh^2(t') \sum_{i+j=s} n_i(t') n_j(t') \, dt'.
\end{equation*}
Equations for the total density and monomers show that the total aggregate concentration $n(t)$ and the monomer concentration $n_1(t)$ relax to a steady state at an exponential rate
$$
1 - n(t) = 1 - \tanh(t) \simeq 2e^{-2t}.
$$
The same holds for other concentrations, but extracting precise analytical expressions for them via recursion is very difficult.

Therefore, it is worthwhile to consider the steady-state problem separately. Assuming $t \rightarrow \infty$, we obtain a system of nonlinear equations:
\begin{equation*}\label{eq:ConstKerStationary}
0 = \frac{1}{2} \sum_{i=1}^{s-1} K n_i n_{s-i} - n_s \sum_{i=1}^{\infty} K n_i + P_s.
\end{equation*}
For the total number of aggregates per unit volume, we have an elementary stationary equation (again for $P=0$, but now with a general constant coagulation kernel $K$):
$$
0 = -\frac{K n^2}{2} + 1,
$$
hence $n = \sqrt{2/K}$. For the generating function $g(z) = \sum_{s=1}^{\infty} z^s n_s$, the stationary equation becomes
$$
0 = \frac{K}{2} g^2 - K g \cdot n + z = \frac{K}{2} g^2 - \sqrt{2K} \, g + z,
$$
from which we take the correct root
$$
g(z) = \frac{\sqrt{2K} \, \bigl(1 - \sqrt{1-z}\bigr)}{K} = \sqrt{\frac{2}{K}} \sum_{s=1}^{\infty} \frac{\Gamma(s-1/2)}{\sqrt{4\pi} \, \Gamma(s+1)} z^s.
$$
Consequently,
$$
n_s = \frac{1}{\sqrt{2\pi K}} \, \frac{\Gamma(s - 1/2)}{\Gamma(s+1)} \approx \frac{1}{\sqrt{2\pi K}} \, s^{-3/2}.
$$

This asymptotic behavior is consistent with the formal accumulation of aggregate mass $\sum_{s \geq 1} s n_s$ in the system over time. Indeed, under zero initial conditions, mass conservation in the presence of a source leads to accumulation: $\sum_{s \geq 1} s n_s(t) = t$ for finite times, so the total mass diverges as $t \to \infty$. 

For large $t$, the concentrations $n_s(t)$ are very close to their steady-state values for $s \ll s_*$, while for $s \gg s_*$, the concentrations $n_s$ are effectively zero. The cutoff mass $s_*$ can be estimated from the relation
\begin{equation*}
\label{cross}
t = \sum_{s=1}^{\infty} s n_s(t) \approx \sum_{s=1}^{s_*} s n_s \sim \sum_{s=1}^{s_*} s^{-1/2} \sim \sqrt{s_*}.
\end{equation*}
Thus, we obtain $s_* \sim t^2$. 

In numerical studies of the finite subsystems of aggregation equations, one can estimate the range of particle sizes for which $n_s$ will be close to stationary, taking into account both the exponential relaxation rate of the solution to the steady state over time and the exponential decay of concentrations for sizes $s > s_*$ for a finite formal integration time of the system of differential equations.

Now let the source power be $P \geq 0$ for particles of size $\hat{s}$, and set $K=1$. For the total particle density $n = \sum_{s=1}^{\infty} n_s$, we obtain the equation
\begin{equation*}\notag
\frac{dn}{dt} = -\frac{n^2}{2} + (1 + P).
\end{equation*}
It is then easy to see for $t \rightarrow \infty$ that
$n = \sqrt{2(P+1)} = \varphi^{1/2}$,
where we denote formally $\varphi = 2(P+1)$. The concentrations $n_s$ can be obtained recursively
\begin{equation*}
\begin{array}{ccc}
n_1 = \varphi^{-1/2} & 
n_2 = \frac{1}{2}\varphi^{-3/2}, & 
n_3 = \frac{1}{2}\varphi^{-5/2}, \\[6pt]
n_4 = \frac{5}{8}\varphi^{-7/2}, & 
n_5 = \frac{7}{8}\varphi^{-9/2}, & 
n_6 = \frac{21}{16}\varphi^{-11/2}, \\[6pt]
n_7 = \frac{33}{16}\varphi^{-13/2}, & 
\ldots & 
n_s = a_s \cdot \varphi^{-s + 1/2}.
\end{array}
\end{equation*}
The coefficients $a_s$ are independent of $P \geq 0$; for $P=0$ and $K=1$, the analytical solution is available from above. Hence, for $s < \hat{s}$, we have
\begin{eqnarray}
\notag
n_s &=& a_s \cdot \varphi^{-s+\frac{1}{2}} =
\underbrace{a_s \cdot 2^{-s+\frac{1}{2}}}_{\text{solution for } P=0} \cdot \; (1+P)^{-s+\frac{1}{2}} \\[4pt]
\notag
&=& \dfrac{1}{\sqrt{2\pi}} \, \dfrac{\Gamma(s - \frac{1}{2})}{\Gamma(s+1)} \cdot (1+P)^{-s+\frac{1}{2}} \approx \dfrac{1}{\sqrt{2\pi}} \, s^{-3/2} (1+P)^{-s+\frac{1}{2}}.
\end{eqnarray}

Hence, the stationary concentrations of the small particles decrease exponentially due to the presence of the second source \cite{matveev2020oscillating}. Next, a sharp jump in the concentrations of particles of size $\hat{s}$ can be obtained directly from the stationary equation:
\begin{eqnarray}
\notag
n_{\hat{s}} &=& \dfrac{\sum_{i=1}^{\hat{s}-1} n_i n_{\hat{s}-i}}{2n} + \dfrac{P}{n}
= \dfrac{\sum_{i=1}^{\hat{s}-1} n_i n_{\hat{s}-i}}{2n} + \dfrac{P}{\sqrt{2(1+P)}} \\[4pt]
\notag
&=& \dfrac{1}{\sqrt{2\pi}} \, \dfrac{\Gamma(\hat{s} - \frac{1}{2})}{\Gamma(\hat{s}+1)} \cdot (1+P)^{-\hat{s}+\frac{1}{2}} + \dfrac{P}{\sqrt{2(1+P)}} \\[4pt]
\notag
&\approx& \dfrac{1}{\sqrt{2\pi}} \, \hat{s}^{-3/2} \cdot (1+P)^{-\hat{s}+\frac{1}{2}} + \dfrac{P}{\sqrt{2(1+P)}}.
\end{eqnarray}
For larger aggregates, we have an exact recurrence relation using the concentrations of smaller particles and the total number of aggregates
\begin{equation*}
n_s = \dfrac{\sum_{i=1}^{s-1} n_i n_{s-i}}{2n} = \dfrac{\sum_{i=1}^{s-1} n_i n_{s-i}}{2\sqrt{2(P+1)}}.
\end{equation*}
Similar results can be obtained for an arbitrary $K>0$ using the generating function method with $g(z) = \sum_{s=1}^{\infty} n_s z^s$. The equation for the generating function becomes
$$
0 = \frac{K}{2} g^2 - K g n + z + P z^{\hat{s}}.
$$
The relevant root is then
\begin{eqnarray*}
g(z) &=& \frac{K n - \sqrt{K^2 n^2 - 2K(z + P z^{\hat{s}})}}{K}
= n \left( 1 - \sqrt{1 - \frac{2}{K n^2}(z + P z^{\hat{s}})} \right) \\[4pt]
&=& \sqrt{\frac{2(1+P)}{K}} \left( 1 - \sqrt{1 - \underbrace{\frac{1}{2(1+P)} \cdot (z + P z^{\hat{s}})}_{\mathbf{z}}} \right).
\end{eqnarray*}
As before, we can expand the Taylor series in the variable $\mathbf{z}$, then expand the brackets in powers of $z$ and readily obtain the values of the concentrations $n_s$ for $s < \hat{s}$:
$$
n_s = \sqrt{\frac{2(P+1)}{K}} \cdot \frac{\Gamma(s-1/2)}{\sqrt{4\pi}\,\Gamma(s+1)} \cdot (P+1)^{-s}
= \sqrt{\frac{1}{2K\pi}} \, \frac{\Gamma(s-1/2)}{\Gamma(s+1)} \cdot (P+1)^{-s+\frac{1}{2}}.
$$

The formal variable $\mathbf{z} = \dfrac{z + P z^{\hat{s}}}{1+P}$ contains both $z$ and $z^{\hat{s}}$. Consequently, explicit formulas for the steady-state concentrations $n_s$ for $s > \hat{s}$ become difficult to obtain because the Taylor series for $g(z)$ mixes terms corresponding to powers $z^{k}$ and $z^{k\hat{s}}$ \cite{matveev2025low}. For large aggregates with size exceeding $\hat{s}$, one can use the recurrence relation obtained above, which yields approximate results with machine precision \cite{matveev2020oscillating}.

A general approach to parameterizing the asymptotic behavior of solutions to Smoluchowski-type equations is known as scaling theory for self-similar solutions \cite{Leyvraz, Krapivsky, Galkin2001}:
$$
n_s(t) \simeq f^{-2}(t) \, \Phi\!\left(\frac{s}{f(t)}\right), \qquad t \gg 1,\; s \gg 1,
$$
where $f(t)$ and $\Phi(\cdot)$ depend primarily on the homogeneity indices of the coagulation kernel and on the initial conditions (typically, $n_s$ are monotonic and smooth).

This formalism allows one to study the properties of solutions even when an analytical solution to the Cauchy problem is unavailable. For example, no exact formulas are yet known for the solutions of aggregation kinetics equations with ballistic and diffusion kernels, which are important for many applications. Nevertheless, based on scaling theory, general asymptotic properties of these solutions have been obtained in \cite{hayakawa1987irreversible, krapivsky2012driven, brilliantov2009model}.

The steady-state solutions obtained above for the aggregation equations with multiple particle sources \cite{matveev2020oscillating, matveev2025low} exhibit a complex nonmonotonic form and approach the scaling steady-state solution $n_s \simeq s^{-3/2}$ only for very large $s \gg 1$. For large values of $P > 1$, the highly oscillatory structure in the solution can persist over thousands or even millions of aggregate sizes, necessitating numerical methods that account for all particle sizes to accurately approximate the concentrations $n_s$.

The basic power-law asymptotics $n_s \simeq s^{-3/2}$ for steady-state concentrations in the problem with a single monomer source and a constant aggregation kernel can be generalized to a wide class of $\alpha$-$\beta$ model kernels and to the diffusion kernel with a monomer source. Specifically, for
\begin{equation*}
\frac{d n_s}{dt} = \frac{1}{2} \sum_{i=1}^{s-1} K_{i,s-i} n_i n_{s-i} - n_s \sum_{i=1}^{\infty} K_{s,i} n_i + \delta_{s,1},
\end{equation*}
with $K_{i,j} = K \cdot (i^{\alpha} j^{\beta} + i^{\beta} j^{\alpha})$, $\alpha, \beta > -1$, $|\alpha - \beta| < 1$, the solution satisfies the relation \cite{hayakawa1987irreversible}
$$
n_s \simeq \sqrt{\frac{(1-4(\alpha-\beta)^2) \cos\bigl(\pi (\alpha-\beta)\bigr)}{4\pi K}} \; s^{\frac{-3+\alpha+\beta}{2}}.
$$
Interestingly, for the diffusion kernel $K_{i,j} = (i/j)^{1/3} + (j/i)^{1/3} + 2$, the solution corresponds to the case $\alpha = -\beta = \frac{1}{3}$ \cite{krapivsky2012driven}.

We are not limited to classical coagulation equations but are also interested in a model of irreversible coagulation with sources and sinks of matter \cite{matveev2020oscillating}:
\begin{equation*}\label{eq:SmolSourcesTrunc}
\frac{d n_s}{dt} = \frac{1}{2} \sum_{i+j=s} K_{i,j} n_i n_j - n_s \sum_{i=1}^{N} K_{s,i} n_i + \underbrace{P_s}_{\text{sources}}, \qquad s = 1,2,\ldots, N.
\end{equation*}
The structure of these equations essentially replicates the classical aggregation model, requiring only that the effect of sources be taken into account when using numerical algorithms and when imposing a finite maximum allowable particle size $N$ (larger particles are considered inactive, e.g., due to gravitational sedimentation). Moreover, the inclusion of sources or terms associated with particle fragmentation events significantly expands the range of potentially interesting problems. In particular, it becomes possible to solve not only direct but also inverse problems, such as reconstructing the source functions $P_s(t)$ from observations of particle concentrations $n_s(t)$ \cite{agoshkov2002solution, penenko2020source, zaks2025fast}.

It is interesting to note that for the $\alpha$-$\beta$ aggregation model with $1 > \alpha = -\beta > 1/2$, the existence of stably time-oscillating solutions \cite{Colm, matveev2017oscillations} has been demonstrated in the irreversible aggregation model with sources and a sink for large particles. Despite the formal existence of a stationary solution, it turns out to be unstable, and in practice the solutions correspond to an attractive limit cycle.

A complete analytical proof of the stability of periodic solutions over time for the coagulation coefficients of the $\alpha$-$\beta$ model has not yet been obtained. However, for a number of simplified models of aggregation-fragmentation processes, the emergence of oscillatory regimes as a result of the Andronov--Hopf bifurcation has been demonstrated \cite{budzinskiy2021hopf, fortin2023stability, pego2020temporal}. Note that calculations for time-oscillating solutions require solving a large number of coupled equations with very small time steps and are computationally intensive \cite{matveev2017oscillations}. Otherwise, the trajectories become unstable and the computational scheme fails.

Unfortunately, the Cauchy problem for aggregation equations with sources can be solved analytically (partially or completely) only in very specific cases (e.g., for a constant kernel with monodisperse or exponential initial conditions \cite{Lushnikov1976, Krapivsky}). Furthermore, the use of scaling theory also has a number of limitations, and its predictions generally require additional verification, especially for relatively small $s$ and $t$ \cite{Leyvraz}.

This situation shows that elaboration of efficient computational algorithms is necessary. Asymptotic and self-similar properties of the solutions $n_s(t)$ for $s \gg 1$ and/or $t \gg 1$ can often be related to the properties of the physical system under study \cite{Leyvraz, kang1986long}, but their derivation and verification typically require numerical solutions of large finite subsystems of kinetic equations, for example,
\begin{equation*}\label{eq:SmolbasicTrunc}
\frac{d n_s}{dt} = \underbrace{\frac{1}{2} \sum_{i+j=s} K_{i,j} n_i n_j}_{\text{formation of particles of size } s} - \underbrace{n_s \sum_{i=1}^{N} K_{s,i} n_i}_{\text{disappearance of particles of size } s}, \qquad s = 1,2,\ldots, N.
\end{equation*}
Such finite systems of differential equations for $N \gg 1$ allow us to approximate the solutions of the original infinite system of kinetic equations for sufficiently long model times \cite{kang1986long, mcleod62, kaushik2023steady, Galkin2001}.

\section{Multiparticle aggregation}
\label{par:Multiparticle Aggregation Equations}

\subsection{Additional notations}

For convenience, we introduce some notation. In the vector notations, the discrete coagulation equations correspond to a system of ODEs with nonlinear operator $\textbf{S}$ corresponding to the right-hand side
\begin{equation}\notag
\left\lbrace
\begin{matrix}
\dfrac{d \textbf{n}}{d t} = \textbf{S}(\textbf{n}(t)),\\[6pt]
\textbf{n}(0) = \textbf{n}_0,
\end{matrix}\right.
\end{equation}
where
\begin{eqnarray}\notag
\textbf{n}^{\top}(t) = \begin{bmatrix}
n_1(t) & n_2(t) & \ldots & n_M(t) & \ldots
\end{bmatrix}, \quad 
\textbf{n}_0^{\top} = \begin{bmatrix}
n_{1_0} & n_{2_0} & \ldots & n_{M_0} & \ldots
\end{bmatrix}.
\end{eqnarray}
If particle aggregation due to three-particle collisions is possible, the aggregation equations take the form (see e.g. \cite{stefonishin2018tensor})
\begin{gather*}
\notag
\frac{d n_k}{dt} = S^{(2)}_k(\mathbf{n}) + S^{(3)}_k(\mathbf{n}), \qquad k = 1, 2, \ldots, \infty,\\[4pt]
\notag
S^{(2)}_k(\mathbf{n}) = \frac{1}{2} \sum_{i_1 + i_2 = k} C^{(2)}_{i_1, i_2} n_{i_1} n_{i_2} - n_k \sum_{i_1 = 1}^{\infty} C^{(2)}_{i_1, k} n_{i_1},\\[4pt]
\notag
S^{(3)}_k(\mathbf{n}) = \frac{1}{6} \sum_{i_1 + i_2 + i_3 = k} C^{(3)}_{i_1, i_2, i_3} n_{i_1} n_{i_2} n_{i_3} - \frac{n_k}{2} \sum_{i_1, i_2 = 1}^{\infty} C^{(3)}_{i_1, i_2, k} n_{i_1} n_{i_2}.
\end{gather*}
It makes sense to distinguish separately two important special cases of this system, describing collisions of only one type, i.e., when the kinetic coefficients responsible for the other type of collisions are set to zero. These systems have the form
\begin{gather*}
\label{s2.1a}
\frac{d n_k}{dt} = S^{(2)}_k(\mathbf{n}), \qquad k = 1, 2, \ldots, \infty;\\[4pt]
\label{s2.1b}
\frac{d n_k}{dt} = S^{(3)}_k(\mathbf{n}), \qquad k = 1, 2, \ldots, \infty.
\end{gather*}
Here equations with $S^{(2)}_k(\mathbf{n})$ are the original two-particle Smoluchowski equations. The three-particle aggregation equations can be easily generalized to the case of multiparticle collisions \cite{jiang1989long, jiang1994scaling, jiang1995large, krapivsky1991aggregation}:
\begin{gather*}
\label{s1.1}
\frac{\mathrm{d} \mathbf{n}}{\mathrm{d}t}
= \sum_{d = 2}^{D} \left[
\mathcal{P}^{(d)}\!\left(\mathbf{n}\right)
+ \mathcal{Q}^{(d)}\!\left(\mathbf{n}\right)
\right],
\end{gather*}
where the operators $\mathcal{P}^{(d)} = \left[ p^{(d)}_{1}, p^{(d)}_{2}, \ldots \right]^{T}$ and $\mathcal{Q}^{(d)} = \left[ q^{(d)}_{1}, q^{(d)}_{2}, \ldots \right]^{T}$ are defined componentwise as
\begin{gather*}
\notag 
p^{(d)}_{k}\!\left(\mathbf{n}\right) 
= \frac{1}{d!}\sum_{I_d = k} 
C^{(d)}_{i_1, i_2, \ldots, i_d} \, n_{i_1} n_{i_2} \ldots n_{i_d}, 
\\ 
\notag
I_d = i_1 + i_2 + \ldots + i_d, 
\quad k = 1, 2, \ldots, \infty, 
\quad d = 2, 3, \ldots, D;
\\[4pt]
\notag 
q^{(d)}_{k}\!\left(\mathbf{n}\right) 
= - \frac{n_{k}}{(d-1)!}\sum_{i_1, i_2, \ldots, i_{d-1} = 1}^{\infty} 
C^{(d)}_{i_1, i_2, \ldots, i_{d-1}, k} \, n_{i_1} n_{i_2} \ldots n_{i_{d-1}},
\\ 
\notag
k = 1, 2, \ldots, \infty, 
\quad d = 2, 3, \ldots, D.
\end{gather*}
The underlying assumption here is that all types of collisions are possible depending on the number $d$ of participating particles for $2 \leq d \leq D$. The number $D$ determines the maximum number of aggregates that can participate in a single merger event. Thus, collisions with $d > D$ are considered unlikely and are not taken into account.

For the Cauchy problem with exclusively three-particle collisions and a constant kernel $C^{(3)}_{i_1, i_2, i_3} \equiv 1$, an explicit solution is known \cite{krapivsky1991aggregation} for the monodisperse initial conditions $n_k(0) = \delta_{1,k}$ (where $\delta_{i,j}$ denotes the Kronecker symbol)
\begin{gather*}
\label{s2.3}
n_{k = 2s-1}(t) = \frac{\Gamma(s - 1/2)}{\Gamma(1/2)\,\Gamma(s)} \cdot n^{3/2}(t) \cdot \bigl(1 - n(t)\bigr)^{s-1},\\[4pt]
n(t) = \sqrt{\frac{3}{3 + 2t}}, \qquad n_{k = 2s}(t) \equiv 0, \qquad s = 1, 2, \ldots, \infty.
\end{gather*}
The vanishing of all even-numbered concentrations $n_{k=2s}$ follows naturally from the choice of initial conditions, the absence of fragmentation processes, and the absence of two-particle collisions. It suffices to note that initially all particles are of odd unit size, and in triple collisions of odd-sized particles, the formation of even-sized aggregates is impossible \cite{krapivsky1991aggregation}.

\begin{gather*}
\notag
n_{k = 1 + s(D-1)}(t)
= \frac{\Gamma\!\bigl(k (D-1)^{-1}\bigr)}{\Gamma\!\bigl((D-1)^{-1}\bigr) \, \Gamma(s+1)} 
\cdot \bigl(n(t)\bigr)^{D/(D-1)} \cdot \bigl(1 - n(t)\bigr)^{s}, 
s = 0, 1, 2, \ldots, \infty;
\\[4pt]
\label{s4.4}
n_{k \neq 1 + s(D-1)}(t) \equiv 0, ~ s = 0, 1, 2, \ldots, \infty; ~
n(t) = \left( 1 + \frac{(D-1)^2}{D!} \cdot t \right)^{-1/(D-1)}.
\end{gather*}
The vanishing of all concentrations $n_{k \neq 1 + s(D-1)}$ again follows naturally from the choice of initial conditions, the absence of fragmentation processes, and the presence of only $D$-particle collisions.

\subsection{Source-induced three-particle problem}

In the case of the monomer source problem for strictly ternary aggregation ($K_{i,j} \equiv 0$) with constant fusion rates, we set $K_{i,j,k} = 6$ and $J = 2$, rescaling the concentration and time units accordingly \cite{krapivsky2024gelation}. The kinetic equations then take the form
\begin{equation*}
\label{ckt-3}
\frac{d n_s}{dt} = \sum_{i+j+k=s} n_i n_j n_k - 3 n_s n^2 + 2\delta_{s,1}.
\end{equation*}
The total cluster concentration evolves according to
\begin{equation*}
\label{ct-3}
\frac{d n}{dt} = -2n^3 + 2.
\end{equation*}
The steady-state value $n = 1$ explains the specific choice of the monomer source power. An exact (albeit implicit) time-dependent solution of this equation, given an initially empty system $n(0)=0$, can in principle be found:
\begin{equation*}
\label{ct-3-imp}
\ln\frac{\sqrt{1+n+n^2}}{1-n} + \sqrt{3}\tan^{-1}\!\left(\frac{1+2n}{\sqrt{3}}\right) = 6t + \frac{\pi\sqrt{3}}{6}.
\end{equation*}

Next, we may consider the case where both binary and ternary aggregation events with constant kernels are allowed simultaneously, following \cite{krapivsky2024gelation}. The governing kinetic equations can then be written as
\begin{eqnarray*}
\label{ckt-23}
\frac{d n_s}{dt} &=& \lambda \sum_{a+b=s}n_a n_b-2\lambda n_s n \nonumber\\
&+& \sum_{i+j+k=s}n_i n_j n_k -3n_s n^2+(2+\lambda)\delta_{s,1}.
\end{eqnarray*}

Without loss of generality, we set the rate constant for pairwise aggregation events to $\lambda>0$, while the intensity of triple fusion events is set to unity. For convenience, the monomer source term is chosen such that the total steady-state concentration of particles equals unity ($n(\infty)=1$), as can be verified from the corresponding equation
\begin{equation*}
\label{ct-3lam}
\frac{d n}{dt}=-\lambda n^2 - 2n^3+(2+\lambda).
\end{equation*}

At steady state, the original system of differential equations reduces to
\begin{equation*}
\label{ck-23}
\lambda \sum_{a+b=s}n_a n_b+\sum_{i+j+k=s}n_i n_j n_k +(2+\lambda)\delta_{s,1}=(2\lambda+3)n_s.
\end{equation*}

Exact expressions can be obtained for the first several steady-state concentrations \cite{krapivsky2024gelation} but analytical formulas quickly become intractable.

The generating function $g(z)$ satisfies the cubic equation
\begin{equation*}
\label{GF:23}
(2\lambda+3) g =\lambda g^2+ g^3 + (2+\lambda)z.
\end{equation*}

Let us introduce the notations
$$\sigma(z) = (-27 z \lambda - 54 z - 2 \lambda^{3} + 18 \lambda^{2} + 27 \lambda)^{2}
+ 4 ( -\lambda^{2} + 6 \lambda + 9)^{3},$$
$$\omega(z) = - 27 z \lambda - 54 z - 2 \lambda^{3} + 18 \lambda^{2} + 27\lambda,$$
$$\gamma(z) = \frac{\sqrt{\sigma(z)}+\omega(z)}{2}.$$
The unique real root of the cubic then takes the form (obtainable with any computer algebra system)
\begin{gather*}
g(z) = \frac{1}{3} \left( \gamma^{1/3}(z)
-
\frac{-\lambda^{2} + 6 \lambda + 9}{\gamma^{1/3}(z)}
- \lambda \right).
\end{gather*}

Deriving readable series expansion coefficients for $g(z)$ from this solution again appears impractical. Nevertheless, the behavior of the solution for large $s \gg 1$ can be readily deduced from the asymptotic expansion of $g(z)$ near $z=1$ because the equation for generating function yields
\begin{equation*}
1-g(z) \simeq \sqrt{\frac{2+\lambda}{3+\lambda}}\,\, \sqrt{1-z},
\end{equation*}
which leads to the interesting observation
\begin{equation*}
\label{ck-23_final}
n_s \simeq \sqrt{\frac{2+\lambda}{4\pi(3+\lambda)}}\,\, s^{-3/2}
\end{equation*}
for $s\gg 1$. Thus, for aggregation with a source and triple particle fusion events at constant aggregation kernels, the same asymptotic $s^{-3/2}$ behavior is preserved as in the purely binary case (but the prefactor changes). The relaxation to the steady state is again exponential $1 - n(t) \sim A e^{-\lambda t}$ (the exact expression for $A$ and $\lambda$ follows from linearization around $n=1$).

Finding the asymptotic behavior of the concentrations is more challenging. This structure is shown in \cite{krapivsky2024gelation}, even though the concentrations of light clusters exhibit significant oscillations for $\lambda \ll 1$. The binary aggregation process dominates for large aggregates ($s \gg 1$), where the oscillations disappear, the solution decays monotonically in accordance with theory, and the asymptotic behavior $n_s \sim s^{-3/2}$ emerges.

The estimate $s_* \sim t^2$ for the cutoff mass (beyond which concentrations decrease exponentially) justifies the numerical integration of the finite subsystem of differential equations (see \cite{krapivsky2024gelation}). The steady-state solutions exhibit the same power-law decay $n_s \sim s^{-3/2}$ again.

This algebraic decay of concentrations is a fairly general phenomenon: for multiparticle aggregation with constant coefficients and a monomer source, the high-mass tail is characterized by a universal exponent of $3/2$, and only the amplitude and the behavior at low masses depend on the details of the reaction rates. The proof of these statements is a straightforward extension of the generating function analysis described above: one derives the singularity $\sqrt{1-z}$ of the generating function $g(z)$ near $z=1$, which implies $n_s \sim s^{-3/2}$, and then fixes the concentration amplitude (see \cite{krapivsky2024gelation} for detail). 

\section{Aggregation-fragmentation equations}

\subsection{Collisional shattering}

In planetary rings, the formation of new aggregates occurs through pairwise collisions of particles, followed by either their fusion or fragmentation into smaller particles \cite{Cuzzi, Esposito}. The outcome of a collision is determined by the magnitude of the kinetic energy of the relative motion of the colliding particles: at low energy values, the particles merge, while at high energy values, they break into fragments. It can be shown that, over a wide range of parameters, the frequency of collisions leading to particle aggregation is proportional to the frequency of collisions leading to particle fragmentation. Consequently, the kinetic coefficients for fusion $K_{ij}$ and for impact fragmentation $A_{ij}$ differ by a multiplicative factor $\lambda > 0$, i.e., $A_{ij} = \lambda K_{ij}$. Thus, we arrive at the following system of equations describing the dynamics of particle aggregation and fragmentation in planetary rings \cite{PNAS}

\begin{align*}
\label{eq:AggShatt}
\frac{d n_1}{dt} &= \underbrace{- n_1 \sum_{i=1}^{\infty} K_{1,i} n_i}_{\text{aggregation}} + \underbrace{\frac{\lambda}{2}\sum_{i \geq 2}\sum_{j \geq 2} (i+j) K_{i,j} n_i n_j + \lambda \sum_{j \geq 2} j K_{1,j} n_j}_{\text{monomer influx after fragmentation}}, \\[4pt]
\frac{d n_s}{dt} &= \underbrace{\frac{1}{2} \sum_{i+j=s} K_{i,j} n_i n_j - n_s \sum_{i=1}^{\infty} K_{s,i} n_i}_{\text{aggregation}} - \underbrace{\lambda n_s \sum_{i=1}^{\infty} K_{s,i} n_i}_{\text{shattering into monomers}}, \qquad s = 2, 3, \ldots, N, \ldots
\end{align*}
As a characteristic of the state of the system at time $t$, we again use the total aggregate concentration:
$$
n(t) = \sum_{j=1}^{\infty} n_j(t).
$$
In the case $K_{i,j} \equiv 1$ and $n_k(0) = \delta_{1k}$, an exact expression for the dynamics of the total number of aggregates per unit volume is available \cite{PNAS}
\begin{equation*}
\label{eq:norm}
n(t) = \left\{
\begin{array}{ll}
\dfrac{2\lambda}{1 + 2\lambda - \exp(-\lambda t)}, & \text{if } \lambda > 0, \\[10pt]
\dfrac{2}{2 + t}, & \text{if } \lambda = 0.
\end{array}
\right.
\end{equation*}
The evolution of the monomer concentration is given by the less elegant but simple equation
$$
n_1(t) = \frac{\lambda}{1+\lambda} \left(1 + \frac{1}{\lambda} \left( \frac{1+2\lambda}{2\lambda} \cdot e^{\lambda t} - \frac{1}{2\lambda} \right)^{-(2+2\lambda)/(1+2\lambda)} \right).
$$

Analytical expressions for the concentrations $n_s(t)$ for $s \geq 2$ remain unknown for $\lambda > 0$ and can only be obtained approximately using numerical methods \cite{matveev2018anderson, stadnichuk2015smoluchowski}. At the same time, the stationary solutions of these equations have been found in \cite{PNAS}, again using the generating function method. For the generating function $g(z) = \sum_{s=1}^{\infty} n_s z^s$, the equation takes the form
$$
\frac{g^2}{2} - \frac{2\lambda(1+\lambda)}{1+2\lambda} \, g + \frac{2\lambda^2}{1+2\lambda} \, z = 0,
$$
and the root is expressed as
$$
g(z) = \frac{2\lambda(1+\lambda)}{1+2\lambda} \left(1 - \sqrt{1 - \frac{1+2\lambda}{(1+\lambda)^2} \, z} \right),
$$
from which we obtain
$$
n_s = \frac{\lambda}{\sqrt{4\pi}} \, \frac{2(1+2\lambda)^{s-1}}{(1+\lambda)^{2s-1}} \cdot \frac{\Gamma(s-1/2)}{\Gamma(s+1)}.
$$
For $s \gg 1$, we obtain the asymptotic behavior
$$
n_s \simeq \frac{\lambda}{\sqrt{4\pi}} \, s^{-3/2} e^{-\lambda^2 s}.
$$
Thus, for $0 < \lambda \ll 1$ and $s < \lambda^{-2}$, the power-law asymptotics dominates over the exponential decay. This structure of the steady-state solution can be generalized to the case of a more complex kernel \cite{PNAS}
$$
K_{i,j} = K_0 \cdot (i j)^{\mu}, \qquad \mu < 1/2,
$$
where the algebraic part and prefactor change,
$$
n_s \simeq \frac{n(\infty)}{\sqrt{4\pi}} \, e^{-\lambda^2 s} \, s^{-3/2 - \mu}.
$$
Here $n(\infty)$ corresponds to the value of the steady-state total aggregate density. This analytical expression for the solution is not complete, but it is very convenient for testing various numerical procedures used to study steady-state particle size distributions with more complex aggregation-fragmentation kernels. In the case of the $\alpha$-$\beta$ kernel, stable dynamic oscillatory solutions that preserve the mass conservation law can be found for this model \cite{matveev2017oscillations}. Despite strong numerical evidence, a rigorous analytical justification of this effect is not known to date.

\subsection{Aggregation with spontaneous fragmentation }

Aggregation equations that account for spontaneous particle disintegration play an important role in many applied problems \cite{brilliantov2009model, PNAS}:
\begin{equation*}
\label{eq:AggSpontFrag}
\frac{d n_s}{dt} = \underbrace{\frac{1}{2} \sum_{i+j=s} K_{i,j} n_i n_j - n_s \sum_{i=1}^{M} K_{s,i} n_i}_{\text{aggregation}} - \underbrace{F_s n_s + \sum_{k>s} \varphi(s,k) F_k n_k}_{\text{spontaneous fragmentation}},
\end{equation*}
where the coefficients $F_s \geq 0$ specify the intensities of spontaneous breakage events for particles of size $s > 1$, and the values $\varphi(s,i) \geq 0$ correspond to the fractions of fragments of size $s$ produced during the decay of larger particles of size $i > s$. From physical considerations, the simplest properties of the coefficients $\varphi(s,i)$ are clear:
\begin{enumerate}
\item $\varphi(s,k) \equiv 0$ if $s > k$, i.e., the formation of fragments larger than the original particle is impossible;
\item $\sum\limits_{s=1}^{k-1} s \, \varphi(s,k) = k$, i.e., the total mass of fragments from the breakage of a particle consisting of $k$ monomers equals $k$.
\end{enumerate}

Several interpretable parameterizations of $\varphi(s,k)$ are quite easy for observation. In particular, the following variants are of interest:
\begin{enumerate}
\item $\varphi(s,k) = \dfrac{k}{s\cdot(k-1)} $ -- uniform distribution of masses for fragments of sizes $s \leqslant  k-1$;
\item $\varphi(s,k) = \dfrac{2}{k-1}$ -- uniform distribution of the number of debris for $s \leqslant  k-1$;
\item "Half-splitting" fragmentation model:
$$
\varphi(s,k) = \left\lbrace \begin{matrix}
2, & s = \ell, & k = 2\ell, & \text{(even-sized particle)}, \\[4pt]
1, & s = \ell,\, \ell+1, & k = 2\ell+1, & \text{(odd-sized particle)}, \\[4pt]
0, & \text{otherwise},
\end{matrix} \right.
$$
\item $\varphi(s,k) = k \, \delta_{s,1}$ -- spontaneous shattering exclusively into monomers (yielding exactly $k$ monomers).
\end{enumerate}

It is interesting to note that in the first two cases, the coefficients correspond to a separable function for which the separable representation $\varphi(s,k) = \psi(s) \cdot \theta(k)$ is possible, whereas in the third and fourth cases, the matrix of values $\varphi(s,k)$ is very sparse. 

The presence of a sparse, separable, or low-rank structure (i.e., $\varphi(s,k) = \sum\limits_{\alpha=1}^r \psi_\alpha(s) \cdot \theta_\alpha(k)$) in the fragmentation coefficients enables the construction of exceptionally efficient algorithms for computing the spontaneous fragmentation terms in a number of steps linear in the number of kinetic equations \cite{zaks2025fast}. These properties were used in a numerical study of solutions to aggregation-fragmentation equations for soil microparticles \cite{vasilyeva2017microbially} and numerical method for the inverse problem of the source function identification in the coagulation-fragmentation equations \cite{zaks2025fast}.

For the case of a model with spontaneous particle disintegration into monomers, size-independent aggregation rates $K$, fragmentation rates $F$, and monodisperse initial conditions, the generating function method can again be used to obtain a steady-state solution analytically. Indeed, the aggregation-fragmentation equations take the form
\begin{align*}
\frac{d n_s}{dt} &= \frac{1}{2} \sum_{i+j=s} K n_i n_j - n_s \sum_{i=1}^{\infty} K n_i - F n_s, \\
\frac{d n_1}{dt} &= -n_1 \sum_{i=1}^{\infty} K n_i + \sum_{i=2}^{\infty} F i n_i.
\end{align*}
Accounting that $n(t) = \sum\limits_{i=1}^{\infty} n_i$ and $m(t) = \sum\limits_{i=1}^{\infty} i n_i = 1$, we obtain
\begin{align*}
\frac{d n_s}{dt} &= \frac{K}{2} \sum_{i+j=s} n_i n_j - K n_s n - F n_s, \\
\frac{d n_1}{dt} &= -K n_1 n + F \left( \sum_{i=1}^{\infty} i n_i - n_1 \right) = -K n_1 n + F (1 - n_1).
\end{align*}
From this we derive the equation for the total concentration
$$
\frac{d n}{dt} = -\frac{K}{2} n^2 - F n + F,
$$
whose stationary solution is the positive root of the quadratic equation
$$
n(\infty) = \frac{F}{K} \left( \sqrt{1 + \frac{2K}{F}} - 1 \right).
$$
Next, we write the equation for the generating function $g(z) = \sum_{s=1}^{\infty} n_s(t) z^s$
\begin{align*}
\begin{cases}
\dfrac{\partial g(z,t)}{\partial t} = \dfrac{K g^2(z,t)}{2} - K g(z,t) n(t) - F g(z,t) + z F, \\[8pt]
g(z,t=0) = z.
\end{cases}
\end{align*}
From this we obtain the stationary equation
$$
\frac{K g^2}{2} - g \bigl( K n(\infty) + F \bigr) + z F = 0.
$$
Considering the structure of the series for $g(z,t)$, the physically admissible root of this quadratic equation is
$$
g(z) = \frac{\sqrt{F^2 + 2KF}}{K} \left( 1 - \sqrt{1 - \frac{2K}{2K + F} \cdot z} \right).
$$
Expanding in a Taylor series, we obtain the stationary solution
$$
n_s = \frac{\sqrt{F^2 + 2KF}}{K} \cdot \frac{\Gamma(s-1/2)}{\sqrt{4\pi}\,\Gamma(s+1)} \cdot \left( \frac{2K}{2K+F} \right)^s \simeq \sqrt{\frac{F^2 + 2KF}{4\pi K^2}} \; s^{-3/2} e^{-A s},
$$
where $A = \ln\left(\dfrac{2K+F}{2K}\right) > 0$. 

Thus, the $\varepsilon$-support of the solution to the aggregation-fragmentation equation is bounded for any $K, F > 0$ and is determined solely by the relationship between the aggregation and fragmentation rates. If the particle fragmentation rate is small relative to the aggregation rate $0 < \frac{F}{2K} \ll 1$, then obtaining a solution with acceptable accuracy within the formally finite but very large $\varepsilon$-support of the solution will require the use of tens of thousands or even millions of equations. 

The existence and uniqueness of a solution to the Cauchy problem for more general coagulation-fragmentation equations in continuous form,
\begin{equation*}
\begin{aligned}
\frac{\partial n(v,t)}{\partial t} = &\ \frac{1}{2} \int_{0}^{v} K(v-u, u) \, n(v-u,t) \, n(u,t) \, du - n(v,t) \int_{0}^{\infty} K(v,u) \, n(u,t) \, du \\
& + \int_{v}^{\infty} F(u,v) \, n(u,t) \, du - \frac{n(v,t)}{v} \int_{0}^{v} u \, F(v,u) \, du,
\end{aligned}
\label{eq:cont_coag_frag}
\end{equation*}
were studied in the pioneering theoretical works of Melzak \cite{Melzak1957} and McLeod \cite{mcleod64}. These equations have a number of important applications for modeling physical processes \cite{PNAS, Cuzzi, Esposito}, such as the growth of prion particles \cite{ouellet2025mechanical, Prion}, the population dynamics of plankton colonies \cite{slomka2020bursts}, and the dynamics of dispersed aerosol particles in the atmosphere \cite{aloyan97}. In these equations, the function $F(u,v)$ is called the fragmentation kernel and corresponds to the fraction of fragments of size $v$ produced during the fragmentation of a larger particle of size $u \geq v$. The factor $1/v$ in the decay term corresponds to normalization to the total sum of the fragment masses.

\subsection{Low-rank matrices for the truncated aggregation equations}

Following the basic idea from \cite{matveev2014fast, matveev2015fast} we assume that the the aggregation kernel admits an exact or approximate decomposition as matrix
\begin{equation*}
    K_{i,j} = \sum_{\alpha=1}^{r_a} u_{\alpha}(i) \, v_{\alpha}(j), \qquad i,j = 1,\dots,N.
\end{equation*}
with rank $r_a \ll N$. This decomposition can be obtained numerically via truncated singular value decomposition (SVD) or other methods (e.g. randomized SVD of adaptive cross approximation) \cite{tyrtyshnikov1997brief}. 

The term corresponding to collisional decay of concentrations is
\begin{equation*}
    L_s = n_s \sum_{j=1}^{M} K_{s,j} n_j.
\end{equation*}
Substituting the skeleton decomposition we obtain that
\begin{equation*}
    L_s = n_s \sum_{\alpha=1}^{r_a} u_{\alpha}(s) \left( \sum_{j=1}^{N} v_{\alpha}(j) n_j \right).
\end{equation*}
Further, we need to save the $r_a$ sums $S_{\alpha} = \sum_{j=1}^{N} v_{\alpha}(j) n_j$,  for $\alpha = 1,\dots,r_a$ within $O(N r_a)$ operations. Then we finalize the calculation
$$
    L_s = n_s \sum_{\alpha=1}^{r_a} u_{\alpha}(s) S_{\alpha}
$$
that takes another $O(r_a N)$ operations. The gain term reads
$$
    G_s = \frac{1}{2} \sum_{i=1}^{s-1} K_{i,s-i} \, n_i \, n_{s-i}.
$$
Using the low-rank represenration of the aggregation kernel we get
\begin{equation*}
    G_s = \frac{1}{2} \sum_{\alpha=1}^{r_a} 
    \sum_{i=1}^{s-1} 
    \bigl( u_{\alpha}(i) n_i \bigr) \cdot 
    \bigl( v_{\alpha}(s-i) n_{s-i} \bigr).
\end{equation*}
Each inner sum over $i$ is a linear convolution of two sequences. Define
$a_{\alpha}(i) = u_{\alpha}(i) n_i$ and $b_{\alpha}(i) = v_{\alpha}(i) n_i$.
Then
$$
    G_s = \frac{1}{2} \sum_{\alpha=1}^{r_a} \bigl( a_{\alpha} * b_{\alpha} \bigr)_s,
$$
where $(a * b)_s = \sum_{i=1}^{s-1} a_i b_{k-i}$ is the discrete linear convolution. All convolutions can be computed efficiently using the fast Fourier transform (FFT) after zero‑padding the sequences to length $2M$ that introduce necessary symmetry each periodic convolution costs $O(N \log N)$ operations \cite{matveev2014fast, matveev2015fast}. Hence the total cost for all $r_a$ convolutions is $O(r_a \, N \log N)$ as well as final complexity for the calculation of the right-hand side of the truncated system with $N$ equations. Given the solution $n_s^{(k)}$ at time $t^{(k)}$, the forward Euler step reads
\begin{equation*}
    n_s^{(k+1)} = n_s^{(k)} + \Delta t \left( G_s^{(k)} - L_s^{(k)} \right),
    \label{eq:explicit_euler}
\end{equation*}
where the gain $G_s^{(k)}$ and loss $L_s^{(k)}$ are evaluated using the current state $n^{(k)}$  with use of the low-rank expressions from above. Hence, the fully explicit update of concentrations takes just $O(N r_a \log N)$ operations making the application of basic explicit time-integration methods (e.g. any Runge-Kutta rule) with adaptive time-steps schemes for modelling the aggregation processes very attractive \cite{matveev2024adaptive}.

These ideas have been successfully generalized also for the case of the multiparticle  \cite{stefonishin2018tensor, stefonishin2019tensors} and the multicomponent coagulation \cite{matveev2016tensor} via exploitation of the low-rank tensor decompositions such as canonical polyadic \cite{kolda2009tensor} and tensor train formats \cite{oseledets2011tensor}.

\begin{figure}
\includegraphics[width=0.49\textwidth]{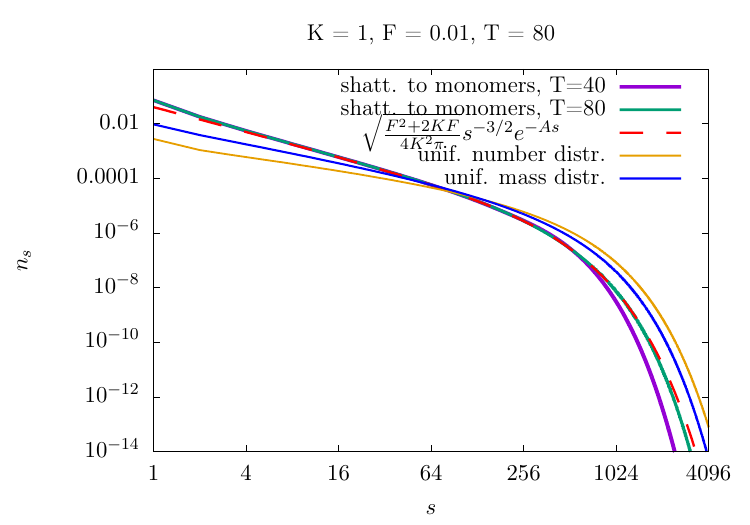}~
\includegraphics[width=0.49\textwidth]{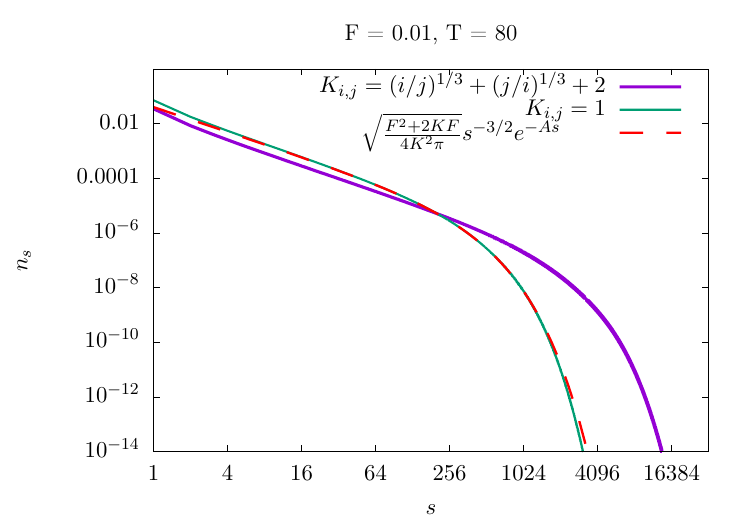}
\caption{Left panel: constant aggregation kernel and size-independent spontaneous fragmentation rates for three models of fragmentation. Right panel: solutions with constant and diffusion aggregation kernels with spontaneous size-independent shattering of into monomers. }
\label{Fig:numerical_results}
\end{figure}

\subsection{Efficient calculation of the spontaneous fragmentation}

Application of matrix methods for fast computations is much less discussed in the case of spontaneous fragmentation. In this section,  we fix this issue and show that for a wide class of the spontaneous fragmentation models this part can be also calculated efficiently.

If one uses $N$ equations in the truncated aggregation fragmentation equations then evaluation of spontaneous fragmenration terms $F_s \cdot n_s$ for each $s = 1,2, \ldots, N$ takes just $N$ operations. However, a straight-forward calculation of the truncated sums 
$$d_s = \sum\limits_{k>s} \varphi(s,k) F_k n_k = \sum\limits_{k=s+1}^{N} \varphi(s,k) F_k n_k$$
in general requires $O(N^2)$ operations and slows-down the whole numerical time-integration process. 

One may assume that symmetrized debris distribution function 
$$\widehat{\varphi}(s,k) = 
\left\lbrace 
\begin{matrix}
\varphi(s,k), & s < k, \\
\varphi(k,s), & s \geq k
\end{matrix} \right., $$
is also low-rank and has a simple skeleton decomposition 
$$\varphi(s,k) = \sum_{\alpha=1}^{r_f} p_{\alpha}(s) q_{\alpha}(k)$$
with $r_f \ll N$. For example, uniform number mass distribution of debris with $\varphi(s,k) = \dfrac{k}{s\cdot(k-1)} $  and their uniform number distribution  $\varphi(s,k) = \dfrac{2}{k-1}$ have exact $r_f=1$. 

Further, we are interested in fast calculation of all vector elements
$\begin{bmatrix}
d_1 & d_2 & \ldots d_s \ldots d_N
\end{bmatrix}^{\top}$ 
from the fragmentation sums. Indeed, in vector-based notations the fragmentation sums receive a simple form
\begin{equation*}
\begin{bmatrix}
d_1 \\ d_2 \\ d_3 \\ \ldots \\  d_s \\ \ldots \\ d_N
\end{bmatrix}
= 
\begin{bmatrix}
0 & \varphi(1,2) F_2 & \varphi(1,3) F_3 & \varphi(1,4) F_4 & \ldots & \varphi(1,N) F_N \\
0 & 0                & \varphi(2,3) F_3 & \varphi(2,4) F_4 & \ldots & \varphi(2,N) F_N  \\
\ldots & \ldots      & \ldots & \ldots & \ldots & \ldots \\
0 & 0 & 0 & 0 & \ldots &  \varphi(N-1, N) F_N  \\
0 & 0 & 0 & 0 & \ldots & 0  \\
\end{bmatrix} \cdot
\begin{bmatrix}
n_1 \\ n_1 \\ n_2 \\ \ldots \\  n_s \\ \ldots \\ n_N
\end{bmatrix}
\end{equation*}
Assume that $r_f=1$ and the symmetrized debris distribution function $\widehat{\varphi}(s,k) = p(s) \cdot q(k)$ is separable, then one may obtain a simple recurrent relation
\begin{align*}
&
d_{N} = 0, \quad d_{N-1} = p(N-1) q(N) F_N n_N = p(N-1) \cdot w_1, \\ 
&
d_{N-2} = p(N-2) q(N-1) F_{N-1} n_{N-1} + p(N-2) q(N) F_N n_N = p(N-2) \cdot (q(N-1) F_{N-1} n_{N-1} + q(N) F_N n_N) \\
& =  p(N-2) \cdot (q(N-1) F_{N-1} n_{N-1} + w_1) = p(N-2) \cdot w_2, \\
&
d_{N-3} = p(N-3) \left( q(N-2) F_{N-2} n_{N-2} + q(N-1) F_{N-1} n_{N-1} + q(N) F_{N} n_{N} \right) \\
& = p(N-2) \cdot \left( q(N-2) F_{N-2} n_{N-2} + w_2 \right) = p(N-2) \cdot w_3, \\
&
\ldots\\
&
d_1 = t_{N-(N-1)} = p(1) \cdot (q(2) F_{2} n_2 + w_{N-1}).
\end{align*}
For arbitrary element in the sequence we should follow the rule
\begin{align*}
d_{N-s} = p(N-s) \cdot (q(N-s+1) F_{N-s+1} n_{N-s+1} + w_{s-1}) = p(N-s) \cdot w_s,\\
\end{align*}
Hence, step-by-step fill of the elements from bottom $t_N$ to the very top $t_1$ requires an elementary update of the cumulant values 
$$w_{s} = w_{s-1} + q(N-s+1) F_{N-s+1} n_{N-s+1}, \quad \text{and } w_0 = 0$$
for each $d_{N-s}$ and additional multiplication by $p(N-s)$. All this work can be done within $O(N)$ operations for whole vector with elements $d_s$, $s = 1, 2, \ldots N$. If the symmetrized debris distribution function $\widehat{\varphi}(s,k)$ matrix is has $r_f$ terms then such trick should be done $r_f$ times and the calculation of the spontaneous fragmentation terms would require just $O(N \cdot r_f)$ operations \cite{tyrtyshnikov2005piecewise}. 

Note, that the values of the total fragmentation rates $F_s$ do not affect the complexity of calculations. Together with the aggregation part, each numerical time-integration step for equations with spontaneous fragmentation events would take just $O(N r_a \log N + N r_f)$ operations for any explicit Runge-Kutta method using $N$ equations (in this work we use the second order Runge-Kutta method with constant time-steps $\Delta t = 10^{-3}$).

\subsection{Numerical experiments}

\begin{figure}
\includegraphics[width=0.49\textwidth]{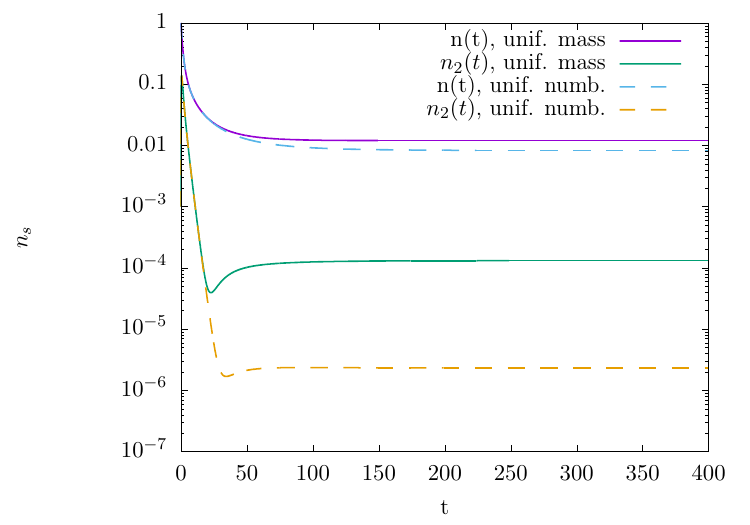}~
\includegraphics[width=0.49\textwidth]{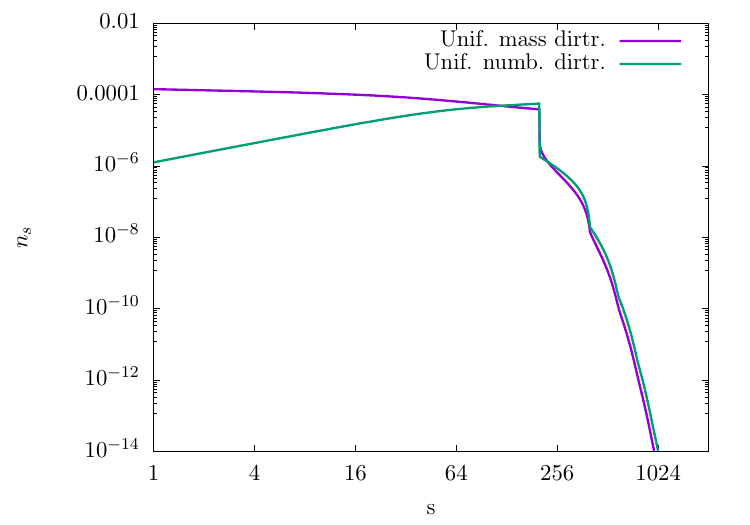}
\caption{Left panel: relaxation dynamics of the concentrations to the steady state for uniform mass and uniform number distributions of debris, size-dependent total fragmentation rates and $K_{i,j}=(i/j)^{0.9} + (j/i)^{0.9}$. Right panel: stationary particle size distributions for these two models of debris distribution. }
\label{Fig:no_oscillatons}
\end{figure}

We verify our analytical expressions for the aggregation-fragmentation equations numerically and present these results in Fig. \ref{Fig:numerical_results}. The variation of models for distribution of debris has a significant influence on the power-law segment (see the left panel of Fig. \ref{Fig:numerical_results} and Fig. \ref{Fig:no_oscillatons}) of the solution but the final exponential cutoff of the concentrations $n_s$ is inevitable. Note, that for the smaller fragmentation rates $F$ the relaxaton time and number of necessary equations for good approximation of the stationary particle size distribution with truncated equations would increase simultaneously: even for relatively large size-independent fragmentation rate $F=0.01$ we see that $T=40$ and $N=2048$ is not sufficient for accurate verification of the analytical solution.

\begin{figure}
\includegraphics[width=0.49\textwidth]{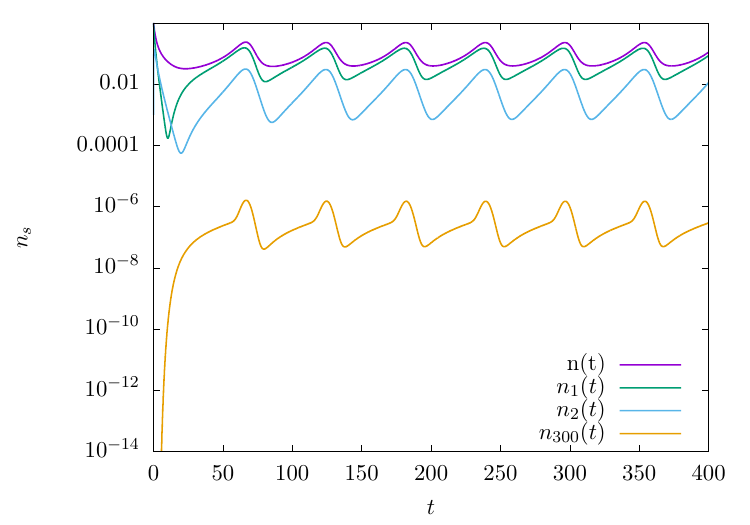}~
\includegraphics[width=0.49\textwidth]{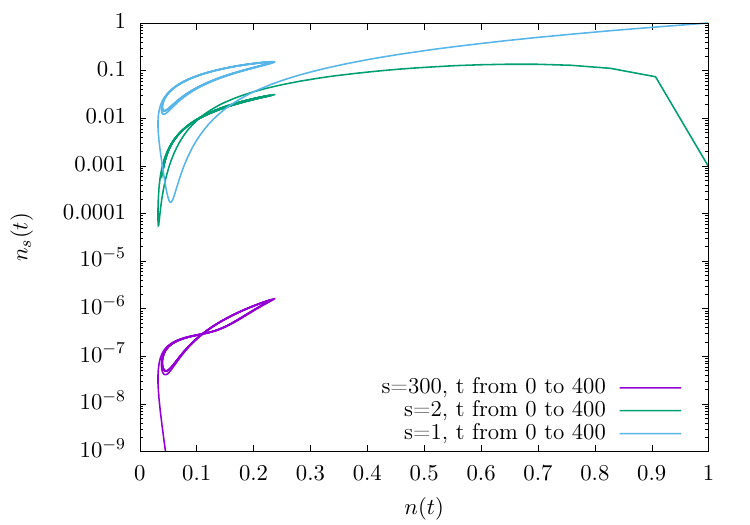}
\caption{Left panel: collective oscillations of particle concentrations with size-dependent spontaneous shattering into monomers and $K_{i,j}=(i/j)^{0.9} + (j/i)^{0.9}$. Right panel: projections of the cyclic phase diagram into coordinates with the total density $n(t)$ as x-axis and $n_{s}(t)$ for $s=1,2, 300$ as y-axis.}
\label{Fig:oscillations}
\end{figure}

On the right panel of Fig. \ref{Fig:numerical_results} we compare the solutions for the constant and diffusion aggregation kernel in case of the constant rates of spontaneous shattering into the monomers. We see the similar algebraic scaling $s^{-3/2}$ of the numerical solutions for these two aggregation kernels with equal homogeneity constant $\rho = 0$. However, more equations are necessary for accurate simulation with the diffusion kernel with more complex structure (we utilized 32768 equations for the diffusion kernel instead of 4096 for the constant kernel) as far as exponential decay of concentrations starts later due to the more intensive aggregation.

The non-decaying temporal oscillations are possible for the aggregation kinetic equations accounting the spontaneous and collisional fragmentation events such as addition-shattering systems \cite{budzinskiy2021hopf, fortin2023stability}. The addition process corresponds to the case of the sparse aggregation coefficients with all zero elements $K_{i,j}$ except for the addition reactions with monomers, where $i = 1$ or $j=1$. More complex analysis for the sparse aggregation coefficients envolving more types of particles is still possible for the temporal oscillations in the Becker-D{\"o}ring equations \cite{pego2020temporal}. In both these cases the numerical time-integration is easy even for the large truncated systems. Linear stability analysis shows that the oscillations arise via the Andronov-Hopf bifurcation as soon as a pair of eigenvalues of the Jacobi matrix at the stationary point receives a positive real part due to variation of parameters in the aggregation-fragmentation rates \cite{budzinskiy2021hopf, fortin2023stability}. Sparsity of the coefficients also makes theoretical treatment of stability of the stationary particle size disributions possible. 

The systems with the fully coupled dense aggregation kernels are tough for theory and require more intensive calculations. Currently, there is a lack of even numerical examples with stable oscillations for the systems with dense aggregation kernels accounting the spontaneous fragmentation terms. In our work we fix this issue taking the $\alpha-\beta$ aggregation kernel with $\alpha=-\beta=0.9$ and following size-dependent total fragmentation rates
\begin{equation*}
F_s = \left\lbrace 
\begin{matrix}
& 0 , & s < 200, \\
& 0.05 \cdot \sqrt{s}, & s \geqslant 200.
\end{matrix} \right.
\end{equation*}

Starting from the monodisperse initial conditions, we observe the non-decaying cyclic collective oscillations of concenrations if the particles shatter only into the monomers (see Fig. \ref{Fig:oscillations}). At the same time for both models of the uniform mass or the uniform number distributions of debris with the same total fragmentation rates we obtain stable relaxation of the numerical solutions to the stationary state (see Fig. \ref{Fig:no_oscillatons}) with rapid exponential cutoff of concentrations $n_s$ of particles of size that exceeds $s=200$, where fragmentation rates become positive and start their growth.

In our paper, as well as in previous works \cite{budzinskiy2021hopf, fortin2023stability, matveev2017oscillations, pego2020temporal}, that also discuss examples with collisional shattering,  the oscillations have been reported for the very narrow debris distributions but for miscellaneous classes of collisional aggregation coefficients that include both sparse and dense kernels. The same fact seems to be true for examples of the source-driven oscillations \cite{Colm}, where the sources are also very localized. Hence, we believe that the oscillations may arise exactly due to this localized structure of sources or debris distributions. Further studies of this phenomena should be directed to more formal approaches allowing to describe necessary or sufficient conditions allowing the birth of the stable temporal oscillations.

\section{Conclusion}

In this work we revisit the generating function approach for obtaining exact analytical solutions of aggregation and fragmentation equations. For discrete coagulation with size independent coefficients under monodisperse initial conditions, we discuss the derivation of explicit closed-form solutions. These results can be extended naturally to continuous formulations, including the two-component case, where exponential initial conditions yield solutions involving modified Bessel functions. We also briefly review cases of additive and product kernels.

Source-driven aggregation produces steady-state distributions characterized by a universal \(s^{-3/2}\) power-law decay and a cutoff mass scaling \(s_{*} \sim t^{2}\). For fragmentation models incorporating spontaneous breakup or collisional shattering, the stationary size distributions become exponentially bounded and can be obtained analytically via the same generating function techniques. Generalizations to three-particle and \(D\)-particle collisions yield solutions expressed in similar terms of gamma functions, with the notable feature that only particles of sizes \(k = 1 + s(D-1)\) are populated.

The generating function method proves to be a powerful analytical framework for studying a wide class of aggregation-fragmentation problems. However, for many practically important kernels — including ballistic, diffusion, and \(\alpha\)-\(\beta\) models — exact analytical solutions remain unknown. In such cases, scaling theory or advanced numerical methods (e.g., low-rank tensor approximations, fast Monte Carlo methods, etc.) are required, particularly for the analysis of dynamic oscillations in aggregation-shattering processes, which still lack rigorous analytical justification.

The analytical results that we revisit here serve as essential benchmarks for validating numerical approaches and provide a basic insight into the universal asymptotic behaviors of aggregating systems, including the persistence of the \(s^{-3/2}\) tail, the exponential cutoff induced by fragmentation, and the oscillatory regimes that challenge conventional steady-state assumptions. 

Finally, we study three models for debris distribution in case of aggregation-fragmentation equations and find the interesting example of temporal oscillations for aggregation-fragmentation equations in case of spontaneous size-dependent shattering of particles into monomers using the efficient numerical method basing on exploitation of the low-rank matrices.

\section*{Acknowledgements}

I am grateful to my senior co-author P. L. Krapivsky who has shown and taught me the techniques of the generating functions during our prior joint publications. I am also grateful to Lutz Warnke for his immediate comments that allowed to extend and polish the initial paper. This work was supported by the Moscow Center of Fundamental and Applied Mathematics at INM RAS (Agreement with the Ministry of Science and Higher Education of the Russian Federation No.075-15-2025-347). 

\bibliographystyle{unsrt}

\bibliography{main}

\end{document}